\documentclass[11pt,a4paper]{article}
\usepackage[left=3cm, right=3cm]{geometry}
\usepackage{tikz}
\usetikzlibrary{trees}
\usepackage[utf8]{inputenc}
\usepackage[T1]{fontenc}
\usepackage[british]{babel}
\usepackage{amsmath}
\usepackage{amsfonts}
\usepackage{amssymb}
\usepackage{makeidx}
\usepackage{graphicx}
\usepackage{color}
\usepackage{enumitem}
\usepackage{authblk}
\usepackage{hyperref}
\usepackage{amsthm}
\usetikzlibrary{backgrounds}
\hypersetup{
	colorlinks=true,
	linkcolor=blue,
	filecolor=blue,      
	urlcolor=blue,
	citecolor=blue
}

\author[1]{Michel Benaïm}
\author[1,2]{Leo Hahn\thanks{Corresponding author:
    {\texttt{leo.hahn@ethz.ch}}.}}
\author[3]{Oliver Tough}
\affil[1]{Institut de mathématiques, Université de Neuchâtel}
\affil[2]{Seminar für Statistik, ETH Zürich}
\affil[3]{Department of Mathematical Sciences, Durham University}
\title{Stochastic extinction: convergence in distribution}

\newtheorem{Lem}{Lemma}

\newtheorem{Def}{Definition}
\newtheorem{Thm}{Theorem}
\newtheorem{Rem}{Remark}
\newtheorem{Prop}{Proposition}
\newtheorem{Asm}{Assumption}

\newcommand{\FM}[1]{d_\mathrm{FM}\left(#1\right)}

\newcounter{mainproofclaim}

\begin{document}
\maketitle

\begin{abstract}
We develop a general framework to show the convergence in distribution of non-irreducible Markov processes to invariant distributions supported on lower-dimensional subsets of the state space. In applications, these singular supports signal the emergence of distinctive long-time behaviours, such as the extinction of species in ecology and the concentration on laminar states in fluid dynamics. The proposed approach is made widely applicable by tailoring its assumptions so that they can be verified using the theory of stochastic persistence/extinction. We illustrate the framework by applying it to stochastic perturbations of the Lorenz and Rosenzweig-MacArthur models and to a general two-dimensional Kolmogorov system.
\end{abstract}

\section{Introduction}

Under, for example, strong Feller and Lyapunov drift conditions, topologically irreducible Markov processes have a unique invariant distribution with full support~\cite{meyn09,benaim22course}. In the absence of irreducibility, however,
there may be multiple invariant measures and their support may be of lower dimension than the state space~\cite{kliemann87,benaim17}. This is crucial in applications because lower-dimensional supports point to the emergence of key, qualitatively distinct phenomena in the long-time regime. In epidemiology, invariant distributions supported on the disease-free set describe the elimination of infection~\cite{DNDY16,benaim19,NYZ20}. In ecology, invariant distributions supported on the boundary of the state space translate the extinction of one or several species~\cite{NY17,HN18,HN18food_chains,HNS22}. In systems of coupled Schrödinger equations~\cite{FH21}, invariant distributions supported on synchronised states describe the emergence of Kuramoto-type collective behaviour. In fluid dynamics, invariant distributions supported on laminar states reflect the decay of turbulence~\cite{CZH21}. Motivated by these applications, this article develops a general approach to show the convergence in distribution of non-irreducible Markov processes to invariant distributions supported on lower-dimensional subsets of their state space. \\

When the process reaches a recurrent class in finite time, the question of long-time behaviour boils down to the irreducible case~\cite{MT93}. On the other hand, when the process remains in a transient class and merely approaches a recurrent class, the tools of stochastic persistence/extinction~\cite{Ben18,HN18,HNC21,FS24+} are required. In that framework, the state space $M$ of a Markov process $(X_t)_{t\ge0}$ is assumed to be the disjoint union of two invariant sets $M_0$ and $M_+$, where $M_0$ is a closed subset of $\partial M_+$. The set $M_0$ typically has lower dimension than $M$, so we think of $M_0$ as translating a specific asymptotic behaviour and of $M_+$ as the generic states. We call $M_0$ the extinction set and $M_+$ the persistence set. Note that, when starting from an initial state in $M_+$, $X_t$ cannot converge to any invariant distribution $\mu_0$ with $\mu_0(M_0)=1$ in total variation. Indeed, the law of $X_t$ and $\mu_0$ remain mutually singular, as $M_0$ and $M_+$ are invariant. One must thus look at weaker forms of convergence~\cite{HMS11,KS18}, which are beyond the classical Harris' ergodic theorem.

For the sake of exposition, first consider the case where $M_0$ carries a single invariant distribution $\mu_0$. In this setting, under suitable additional assumptions, the asymptotic behaviour of $X_t$ is determined by a stochastic version of the Lyapunov exponent: the $H$-exponent $\Lambda$~\cite{BHS08,SBA11}. If $\Lambda > 0$ we speak of persistence. There then exists a unique invariant distribution $\mu_+$ on $M_+$ and, started from any state $x \in M_+$, $X_t$ converges to $\mu_+$ in total variation. If $\Lambda < 0$ we speak of extinction~\cite{FS24+}. In that case, $\mu_0$ is the unique invariant distribution of $(X_t)_{t\ge0}$ and almost surely $d(X_t, M_0) \to 0$ and $\Pi_t \to \mu_0$, where $\Pi_t = \frac1t \int_0^t \delta_{X_s} ds$ is the empirical occupation measure. Crucially, the question of the convergence in distribution of $X_t$ to $\mu_0$ is left open by the general theory. Moreover, the counterexample of Section~\ref{sec:counterexample} shows that $d(X_t, M_0) \to 0$ does not imply the convergence of $X_t$ to $\mu_0$ in distribution, even when $M$ is compact and $\mu_0$ is the sole invariant distribution of $X_t$. \\

In this article, we develop a general framework for upgrading almost sure convergence to the extinction set $M_0$ to convergence in distribution to an invariant distribution supported on $M_0$. When $M_0$ carries a unique, attracting invariant distribution $\mu_0$, it is enough that the process restricted to $M_0$ converges to $\mu_0$ locally uniformly in the initial state. The setting where $M_0$ carries other, non-attracting invariant distributions is handled by requiring uniform convergence only on compact sets that the process visits with high probability. The method can be seen as a generalisation of~\cite{HN18food_chains}. Our contribution is an abstract criterion applicable across different model classes and the convergence to mixtures with initial-state-dependent weights in the case of several attracting invariant distributions on $M_0$. The assumptions are tailored so that they can be verified using the existing theory of stochastic persistence/extinction.

Our framework unifies the model-specific convergence-in-distribution results discussed below, which have been obtained through different approaches. For two-dimensional stochastic population models with quadratic interactions, \cite{RP07,DDT11,DNDY16} use estimates introduced in~\cite{Rud03}. Another method~\cite{NY17,FH21,LG22} constructs an asymptotic coupling between a trajectory in $M_+$ and one in $M_0$. In particular, for the Lotka-Volterra model \cite{NY17}, this yields convergence to a mixture with initial-state-dependent weights. For stochastic food chains, \cite{HN18food_chains,BBN22,Bou23} are similar in spirit to the approach presented in this article. \\

The main result is presented in Section~\ref{sec:main_result}. The counterexample of Section~\ref{sec:counterexample} shows that the almost sure convergence of $d(X_t, M_0) \to 0$ and $\Pi_t \to \mu_0$ does not imply convergence in distribution of $X_t$. Section~\ref{sec:applications} then illustrates the general approach by applying it to 
stochastic perturbations of the Lorenz~\cite{CZH21} and Rosenzweig-MacArthur models~\cite{BCS25} and to a general two-dimensional Kolmogorov system~\cite{HN18}. Finally, Section~\ref{sec:proof_main_result} proves the main result.

\section{Main result}\label{sec:main_result}

We consider a Markov process $(X_t)_{t \ge 0}$ on a Polish space $(M, d)$, sometimes using the notation $(X^x_t)_{t\ge0}$ to emphasise the initial state $x \in M$. Our main result is an abstract framework to upgrade almost sure convergence to the extinction set $M_0$, as provided by stochastic persistence/extinction~\cite{Ben18,HN18,HNC21,FS24+}, to convergence in distribution to an invariant distribution supported on $M_0$. For the sake of exposition, we first state our result in the setting where there is a single invariant distribution $\mu_0$ on $M_0$. This is for example the case for the coupled stochastic Schrödinger equations~\cite{FH21} and the turbulence model~\cite{CZH21}. Essentially, if $d(X_t, M_0) \to 0$ almost surely and the process converges to $\mu_0$ locally uniformly in the initial state on $M_0$, then $X_t$ converges to $\mu_0$ in distribution. The following distance, which metrises convergence in distribution, is key to formulating the uniform-convergence condition. 

\begin{Def}\label{def:FM} The \emph{Fortet-Mourier distance} between two probability measures $\mu$ and $\nu$ is given by
	$$
	d_\mathrm{FM}(\mu, \nu) := \sup \{ \mu(f) - \nu(f) : \Vert f \Vert_{\infty} \le 1 \text{ and } \Vert f \Vert_\mathrm{Lip} \le 1  \}
	$$
	where $\Vert f \Vert_\mathrm{Lip} = \sup_{x \ne y} \frac{|f(x) - f(y)|}{d(x, y)}$.
\end{Def}

Throughout this article, we assume the following slightly strengthened version of Feller continuity.

\begin{Asm}\label{asm:CbM_Feller}
	The map $(t, x) \mapsto P_t f(x) = \mathbb E[f(X^x_t)]$ is jointly continuous for every bounded, continuous function $f\in C_b(M)$.
\end{Asm}

\begin{Thm}\label{thm:single_measure}
	Fix $x \in M$. In addition to the standing Assumption~\ref{asm:CbM_Feller}, assume that:
	\begin{enumerate}[label=(\Alph*\textsubscript{\theThm}), start=2, ref=(\Alph*\textsubscript{\theThm})]
		\item\label{asm:tighness_single_measure} There exists $T > 0$ such that $(X^x_{kT})_{k \ge 1}$ is tight.
		
		\item\label{asm:cv_to_M0_single_measure} There exists a closed set $M_0$ such that
		\begin{equation*}
		d(X_t^{x}, M_0) \to 0 \text{ almost surely}.
		\end{equation*}
		\item\label{asm:cv_to_mu0_single_measure} There exists a probability measure $\mu_0 \in \mathcal P(M_0)$ such that for all compact sets $K \subset M_0$
		\begin{equation}\label{eq:uniform_cv_on_M0}
		\lim_{t \to +\infty} \sup_{y \in K}\FM{\delta_y P_t, \mu_0} = 0.
		\end{equation}
	\end{enumerate}
	Then $X^x_t$ converges in distribution to $\mu_0$.
\end{Thm}

In population dynamics, there is generally more than one invariant distribution on the extinction set. Indeed, in that setting, $(X_t)_{t \ge 0}$ takes its values in $M = \mathbb R_+^n$ and the sets $M_0^I = \bigcap_{i \in I}\{ x_i = 0\}~(I \subset \{ 1, \ldots, n \})$ are all invariant and all possibly carry a different invariant distribution. In particular, $\delta_{(0,  \ldots, 0)}$ is always invariant. This motivates the following result, which includes Theorem~\ref{thm:single_measure} as a special case.

\begin{Thm}\label{thm:main}
Fix $x \in M$. In addition to the standing Assumption~\ref{asm:CbM_Feller}, assume that:
	\begin{enumerate}[label=(\Alph*\textsubscript{\theThm}), start=2, ref=(\Alph*\textsubscript{\theThm})]
		\item\label{asm:tightness} There exists $T > 0$ such that $(X^x_{kT})_{k \ge 1}$ is tight. Let the family of compact sets $(K_\delta)_{\delta > 0}$ be such that
		\begin{align*}
			\inf_{k \ge 1} \mathbb P(X^x_{kT} \in K_\delta) \ge 1 - \delta \text{ for all } \delta> 0.
		\end{align*}
		\item\label{asm:cv_to_M0} The extinction set is a union $M_0 = \cup_{i = 1}^n M_0^i$ of closed sets $M_0^1, \ldots, M_0^n$ such that
		$$
		M_0^i \cap M_0^j \cap K_\delta = \emptyset \text{ for all } i \ne j \text{ and } \delta> 0,
		$$
		and the events $\mathcal A^x_i = \left\{ \lim_{t \to +\infty} d(X^x_t, M_0^i) = 0\right\}$ satisfy
		$$
		\sum_{i = 1}^n \mathbb P(\mathcal A^x_i) = 1.
		$$
		\item\label{asm:cv_to_mu0} There exist probability measures $\mu_0^1 \in \mathcal P(M_0^1), \ldots, \mu_0^n \in \mathcal P(M_0^n)$ such that
		$$
		\lim_{t \to +\infty} \sup_{y \in K_\delta \cap M_0^i} \FM{\delta_y P_t, \mu_0^i}  = 0 \text{ for all } i \in \{1, \ldots, n\} \text{ and }  \delta> 0.
		$$
	\end{enumerate}
	Then $X^x_t$ converges in distribution to $\sum_{i = 1}^n \mathbb P(\mathcal A^x_i) \mu_0^i$.
\end{Thm}

\begin{Rem}\label{rem:main_result}
\begin{enumerate}[label=(\roman*),ref=\theRem (\roman*)]

\item An important special case of Theorem~\ref{thm:main} is when $n = 1$. This is the case where there is a single attracting invariant distribution~$\mu_0 = \mu_0^1$ but possibly also other, non-attracting invariant distributions on the extinction set $M_0 = M_0^1$. Assumption~\ref{asm:cv_to_M0} then coincides with Assumption~\ref{asm:cv_to_M0_single_measure} (almost sure convergence to $M_0$) and Assumption~\ref{asm:cv_to_mu0} becomes
\begin{equation*}
	\lim_{t\to+\infty} \sup_{y \in K_\delta \cap M_0} \FM{\delta_y P_t, \mu_0} = 0.
\end{equation*}
The key difference with Assumption~\ref{asm:cv_to_mu0_single_measure} is that the supremum above is not over an arbitrary compact subset of $M_0$ but over $K_\delta \cap M_0$. In other words, it is enough that the uniform convergence to $\mu_0$ hold on compact subsets of $M_0$ that the process actually visits. This is crucial as Assumption~\ref{asm:cv_to_mu0_single_measure} fails on the support of ergodic measures on $M_0$ that differ from $\mu_0$.

\item Under the assumptions of Theorem~\ref{thm:main}, one has $\mathbb P(\mathcal A_i^x \cap \mathcal A_j^x) = 0$ for $i \ne j$ (see the proof of Claim~\ref{clm:main-convergence}).

\item The convergence to the mixture $\sum_{i=1}^n \mathbb P(\mathcal A_i^x)\mu_0^i$ is analogous to the classical convergence result~\cite[Thm.~8.1]{MT93} when the process enters a recurrent class in finite time. Note that the weights $\mathbb P(\mathcal A_i^x)$ depend on the initial state $x \in M$.
\end{enumerate}
\end{Rem}

The following preliminary definitions are needed to provide a version of Theorem~\ref{thm:main} with conditions that are convenient to verify in practice.

\begin{Def}
	A measurable set $E \subset M$ is \emph{forward invariant} if
	$$
	\mathbb P(X^x_t \in E) = 1 \text{ for all } x \in E \text{ and } t \ge 0.
	$$
\end{Def}

\begin{Def}
	Let $E \subset M$ be forward-invariant. We say that $x \in E$ is \emph{accessible from $E$} if for all open $U \subset M$ containing $x$ and all $y \in E$ there exists $t \ge 0$ such that $\mathbb P(X^y_t \in U) > 0$.
\end{Def}

\begin{Def}
	Let $E \subset M$ be forward-invariant. We say that $x \in E$ is a \emph{Doeblin point for $(X_t)_{t\ge0}$ restricted to $E$} if there exists $m \in \mathcal P(M)$, an open set $U \subset M$ containing $x$ and $a, t	 > 0$ such that
	$$
	\mathbb P( X^y_t \in A) \ge a\,m(A	)
	$$
	for all $y \in U \cap E$ and all measurable $A \subset M$.
\end{Def}

\begin{Thm}\label{thm:practical}
	In addition to the standing Assumption~\ref{asm:CbM_Feller}, assume that:
	\begin{enumerate}[label=(\Alph*\textsubscript{\theThm}), start=2, ref=(\Alph*\textsubscript{\theThm})]
		\item \label{asm:tightness_practical} There exists a forward-invariant $M_{++}$ and a continuous, proper $W : M_{++} \to \mathbb R_+$ such that
		$$
		P_T W(x) \le \rho W(x) + \kappa \text{ for all } x \in M_{++},
		$$
		where $T, \kappa> 0$ and $\rho \in (0, 1)$.
		\item \label{asm:cv_to_M0_practical} The extinction set is a union $M_0 = \cup_{i = 1}^n M_0^i$ of forward-invariant closed sets $M_0^i$ such that
		$$
		M_0^i \cap M_0^j \cap M_{++} = \emptyset \text{ for all } i \ne j,
		$$
		and the events $\mathcal A^x_i = \left\{ \lim_{t \to +\infty} d(X^x_t, M_0^i) = 0\right\}$ satisfy
		$$
		\sum_{i = 1}^n \mathbb P(\mathcal A^x_i) = 1 \text{ for all } x \in M_{++}. 
		$$
		\item \label{asm:cv_to_mu0_practical} For every $i = 1, \ldots, n$, there exists a point $x_i \in M_0^i \cap M_{++}$ that is both accessible from $M_0^i \cap M_{++}$ and a Doeblin point for $(X_t)_{t \ge 0}$ restricted to $M_0^i \cap M_{++}$.
	\end{enumerate}
	The restriction of $(X_t)_{t\ge 0}$ to each $M_0^i \cap M_{++}$ then admits a unique invariant distribution $\mu_0^i  \in \mathcal P(M_0^i \cap M_{++})$ and for all $x \in M_{++}$
	$$
	X^x_t \text{ converges to } \sum_{i = 1}^n \mathbb P(\mathcal A_i^x) \mu_0^i \text{ in distribution}. 
	$$
\end{Thm}

The following remark lays out how Theorem~\ref{thm:practical} fits into the existing theory of stochastic persistence/extinction~\cite{Ben18,HN18,HNC21,FS24+} and how to verify its assumptions. See  also the applications of Section~\ref{sec:applications} for an illustration.

\begin{Rem}\label{rem:how_to_verify_assumptions}
\begin{itemize}
\item By dominated convergence, to verify Assumption~\ref{asm:CbM_Feller}, it is enough that for every $(\mathbb R_+ \times M)$-valued converging sequence $(t_n, x_n) \to (t, x)$, one has $X^{x_n}_{t_n} \xrightarrow{n \to +\infty} X^x_t$ almost surely. Alternatively, one can study the regularity properties of the Kolmogorov backward equation using PDE tools.

\item Applying the theory of stochastic extinction/persistence to the extinction set $M_{00} = M \setminus M_{++}$ provides an efficient way to construct the $W$ of  Assumption~\ref{asm:tightness_practical}, see e.g.~\cite[Lem.~4.12 and Prop.~8.3]{Ben18}.

\item Assumption~\ref{asm:cv_to_M0_practical} can also be shown using stochastic persistence/extinction.

\item The accessibility in Assumption~\ref{asm:cv_to_mu0_practical} typically follows from controllability as in~\cite[Prop.~6.33]{benaim22course}. The Doeblin property typically comes from Hörmander-type conditions as in~\cite[Thm.~6.37]{benaim22course}.
\end{itemize}
\end{Rem}

\section{Instructive counterexample}\label{sec:counterexample}

The following  counterexample shows that, even when $M$ is compact and there is a single invariant distribution $\mu_0$ on the extinction set $M_0$, almost sure convergence to $M_0$ does not imply convergence in distribution to $\mu_0$.
	
\begin{Def}[Counterexample]\label{def:counterexample}
	Consider the state space $M = [0, 1] \times \mathbb S^1$ and the Markov process $X_t = (x(t), \theta(t))$ given by the (deterministic) ODE
	\begin{align*}
		\dot x &= -x^2, & \dot \theta = 1 - \cos \theta + x.
	\end{align*}
\end{Def}

Consider the extinction set $M_0 = \{ 0 \} \times \mathbb S^1$, its single invariant distribution $\mu_0 = \delta_{(0, 0)}$ and the corresponding persistence set $M_+ = M \setminus M_0$.

\begin{Prop}
	Fix $x \in M_+$. One has
	\begin{enumerate}[label=(\roman*)]
		\item $\lim_{t \to +\infty} d(X^x_t, M_0) = 0$ almost surely,
		\item $\Pi^x_t = \frac1t \int_0^t \delta_{X^x_s} ds$ converges weakly to $\delta_{(0, 0)}$, almost surely,
		\item $X^x_t$ does not converge to $\delta_{(0, 0)}$ in distribution.
	\end{enumerate}
\end{Prop}

\begin{proof}
	(i) Solving $\dot x = -x^2$ yields $x(t) = 1/(t+1/x(0))$ and thus $\lim_{t \to +\infty} d(X^x_t, M_0) = 0$.
	
	 (ii) The state space $M$ is compact so $(\Pi_t)_{t\ge0}$ is tight. By~\cite[Thm.~2.2 (ii)]{Ben18} all limit points of $(\Pi_t)_{t\ge 0}$ are invariant distributions of $(X_t)_{t\ge0}$, so $\delta_{(0,0)}$ is the only limit point and hence $\Pi_t$ converges weakly to $\delta_{(0, 0)}$, almost surely.
	
	 (iii) One has
	$$
	\int_0^{+\infty} 1 - \cos \theta(t) + x(t) dt \ge \int_0^{+\infty} \frac{1}{t+ 1/x(0)} dt = +\infty
	$$
	so $\theta(t) = \theta(0) + \int_0^t 1 - \cos \theta(s) + x(s) ds$ goes around $\mathbb S^1$ infinitely many times and does not converge.
\end{proof}

The issue in Counterexample~\ref{def:counterexample} is that, depending on the initial state $x \in M_0$, the convergence of $X^x_t$ to $\mu_0$ can be arbitrarily slow. In other words, Assumption~\ref{asm:cv_to_mu0} is not verified. Similarly, the stable stochastic cycle case in~\cite{BRY25+} shows that $\Pi^x_t$ can fail to converge, even when $d(X^x_t, M_0) \to 0$ almost surely.

\section{Applications}\label{sec:applications}

Before applying Theorem~\ref{thm:practical} to stochastic perturbations of the Lorenz and Rosenzweig-MacArthur models, as well as a stochastic Kolmogorov system, we outline how to use Lyapunov functions to verify Assumption~\ref{asm:tightness_practical}.

\begin{Def}\label{def:extended_generator}
	Let $E \subset M$ be forward-invariant and denote $C(E)$ the set of continuous (possibly unbounded) functions on $E$. Let $(X_t)_{t\ge 0}$ be a Markov process satisfying the standing Assumption~\ref{asm:CbM_Feller}. The restriction $(D(\mathcal L_e^E), \mathcal L_e^E)$ of its \emph{extended generator} to $E$ is the unique linear operator $\mathcal A : C(E) \supset D(\mathcal A) \to C(E)$ such that for all $f, g \in C(E)$
	$$
	f \in D(\mathcal A) \text{ and } \mathcal A(f) = g \iff f(X^x_t) - \int_0^t g(X^x_s) ds \text{ is a local martingale for all } x \in E.
	$$
\end{Def}

When $E$ is clear from context, we simply use the notation $(\mathcal L_e, D(\mathcal L_e))$.

\begin{Def}\label{def:lyapunov} 
Let $E \subset M$ be forward-invariant. We say that $W : E \to \mathbb R_+$ is a \emph{Lyapunov function} for a Markov process with extended generator $(\mathcal L_e^E, D(\mathcal L_e^E))$  if $W \in D(\mathcal L_e^E)$ is proper and
	$$
	\mathcal L_e W \le -a W + C
	$$
	where $a, C > 0$.
\end{Def}

\begin{Prop}\label{prop:decay}
Let $E \subset M$ be forward invariant and let $W : E \to \mathbb R_+$ be a Lyapunov function. For all $x \in E$ and $t\ge 0$ one has
$$
P_t W(x) \le e^{-at} W(x)+ \frac Ca \left(1 - e^{-at}\right).
$$
\end{Prop}

\begin{proof}
Follows by mimicking the proof of~\cite[Thm.~6.1]{meyn93}.
\end{proof}

\subsection{Stochastic perturbation of the Lorenz system}\label{sec:Lorenz}

    Coti Zelati and Hairer~\cite{CZH21} considered the following stochastic perturbation of the Lorenz system, a paradigmatic model for turbulence
	\begin{align}\label{eq:Lorenz_SDE}
			dx = \sigma (y - x)dt, \qquad d y = (x (\rho - z) - y)dt, \qquad d z = (-\beta z + x y)dt  + \hat\alpha dB_t,
	\end{align}
    where $\sigma, \beta, \hat \alpha > 0$ and $\rho \in (0, 1)$. The SDE~\eqref{eq:Lorenz_SDE} takes its values in $M = \mathbb R^3$ and we denote its solution with initial state $(x, y, z) \in M$ by $X^{(x, y, z)}_t = (x(t), y(t), z(t))$. Well-posedness of~\eqref{eq:Lorenz_SDE} and Assumption~\ref{asm:CbM_Feller} follow from classical SDE theory using the Lyapunov function $W : M \to \mathbb R_+$ given by
    $$
        W(x, y, z) = \rho x^2 + \sigma y^2 + \sigma(z - 2\rho)^2.
    $$
    On the invariant subset $M_0 = \{ x = y = 0 \}$, \eqref{eq:Lorenz_SDE} is simply an Ornstein-Uhlenbeck process. The restriction of~\eqref{eq:Lorenz_SDE} to $M_0$ thus admits a unique, Gaussian invariant distribution $\mu_0$. The persistence set $M_+ = M \setminus M_0$ and the extinction set $M_0$ correspond to turbulent and laminar states respectively. Set
    \begin{align*}
    	\chi &=\frac2{1+\sigma}, & \eta &= \frac{1+\sigma}{2\sigma}, & \gamma &= \chi \beta, & \nu^2 &= \chi^5 \sigma, & \alpha &= \nu \sqrt\sigma \hat\alpha,
    \end{align*}
    and $\zeta_\star = 2 + \chi^2 \sigma (\rho - 1)$. Coti Zelati and Hairer~\cite{CZH21} show that the existence of other invariant distributions is governed by the $(\mathbb S^1 \times \mathbb R)$-valued SDE
    \begin{align}\label{eq:Lorenz_aux_SDE}
        d \theta &= (1 - \zeta \sin^2(\theta)) dt, & d \zeta &= -\gamma(\zeta - \zeta_\star) dt + \alpha dB_t.
    \end{align}
	Indeed, if $\mu_\alpha$ denotes the unique invariant distribution of~\eqref{eq:Lorenz_aux_SDE}, then the existence or lack of a second invariant distribution is determined by the sign of the $H$-exponent
    $$
    \Lambda = \int_{\mathbb S^1\times \mathbb R} \left(-1 + \frac \zeta2 \sin(2\theta) \right)d\mu_\alpha(\theta, \zeta).
    $$
	
	\begin{Prop}[\cite{CZH21}]
	\begin{enumerate}[label=(\roman*)]
		\item If $\Lambda < 0$ then $\mu_0$ is the unique invariant distribution of~\eqref{eq:Lorenz_SDE}.
		\item If $\Lambda > 0$ then, in addition to $\mu_0$, \eqref{eq:Lorenz_SDE} admits exactly one other ergodic measure $\mu_+$. This additional ergodic measure satisfies $\mu_+(M_+) = 1$. 
	\end{enumerate}
	\end{Prop}

    Theorem~\ref{thm:practical} yields the following asymptotic behaviour.

	\begin{Prop}
    Let $(x, y, z) \in M_+$.  
		\begin{enumerate}[label=(\roman*)]
            \item If $\Lambda < 0$ then $X^{(x, y, z)}_t$ converges to $\mu_0$ in distribution when $t \to +\infty$.
			\item If $\Lambda > 0$ then $X^{(x, y, z)}_t$ converges to $\mu_+$ in total variation when $t \to+\infty$.
		\end{enumerate}
	\end{Prop}
	
	\begin{proof}
	(i) It suffices to verify the assumptions of Theorem~\ref{thm:practical} with $n = 1$ and $M_{++} = M$. Applying Proposition~\ref{prop:decay} to $W$ yields~\ref{asm:tightness_practical}.	Assumption~\ref{asm:cv_to_M0_practical} is equivalent to $(x(t), y(t)) \to (0, 0)$, which is shown in the proof of~\cite[Thm.~4.1]{CZH21}. 
Finally, the restriction of~\eqref{eq:Lorenz_SDE} to $M_0$ is the Ornstein-Uhlenbeck process, which is elliptic. Hence every point in $M_0$ is accessible from $M_0$~\cite[Prop.~6.33]{benaim22course} and a Doeblin point for $(X_t)_{t\ge 0}$ restricted to $M_0$ by~\cite[Thm.~6.37]{benaim22course}, yielding Assumption~\ref{asm:cv_to_mu0_practical}.
	
    (ii) Every point in $M_+$ is accessible from $M_+$ by~\cite[Prop.~3.3]{CZH21} and~\cite[Prop.~6.32]{benaim22course}, as well as a Doeblin point for $(X_t)_{t \ge 0}$ restricted to $M_+$ by~\cite[Prop.~3.2]{CZH21} and~\cite[Thm.~6.37]{benaim22course}. The proof of~\cite[Thm.~4.1]{CZH21} yields a Lyapunov function for $(X_t)_{t \ge 0}$ restricted to $M_+$, hence the desired convergence in total variation follows from~\cite[Thm.~8.15]{benaim22course}.
\end{proof}
	
\subsection{Rosenzweig-MacArthur prey-predator model}\label{sec:Rosenzweig-MacArthur}

    The Rosenzweig-MacArthur model~\cite{BCS25} is the SDE
    \begin{align}\label{eq:Rosenzweig_MacArthur_SDE}
        dx_1 &= x_1\left(1 - \frac{x_1}\kappa - \frac{x_2}{1 + x_1}\right)dt + \epsilon x_1dB_t, & dx_2 &= x_2\left(-\alpha + \frac{x_1}{1+x_1}\right)dt,
    \end{align}
    where $\kappa, \epsilon, \alpha> 0$. It describes the joint evolution of a population of prey $x_1$ and predators $x_2$. We consider~\eqref{eq:Rosenzweig_MacArthur_SDE} on the state space $M = \mathbb R_+^2$ and denote its solution with initial state $(x_1, x_2)$ by $X^{(x_1, x_2)}_t = (x_1(t), x_2(t))$. Well-posedness of~\eqref{eq:Rosenzweig_MacArthur_SDE} and Assumption~\ref{asm:CbM_Feller} follow from classical SDE theory using the Lyapunov function $U : M \to \mathbb R_+$ given by
    $$
    U(x_1, x_2) = (x_1 + x_2)^4 + 1.
    $$
  The persistence set $M_+ = \mathbb R_{>0}^2$ corresponds to both species surviving, whereas the extinction set $M_0 = M \setminus M_+$ corresponds to at least one species dying out. Importantly, the Brownian noise is only present in the SDE for $x_1$, making~\eqref{eq:Rosenzweig_MacArthur_SDE} degenerate. This models a population of prey that is subjected to random environmental fluctuations and a population of predators that is impacted by the environmental fluctuations only indirectly through its interactions with the prey.
  
  Benaïm, Colombo and Strickler~\cite{BCS25} show that the asymptotic behaviour of $X^{(x_1,x_2)}_t$ is determined by the $H$-exponent
    $$
    \Lambda = \int_{M_0} \left(\frac{x_1}{1 + x_1} - \alpha\right) d\mu_0(x_1, x_2)
    $$
    defined when $\epsilon< \sqrt 2$. Here $\mu_0 = \gamma_{\epsilon, \kappa}(x)dx \otimes \delta_0$ where $\gamma_{\epsilon,\kappa}(x) = x^{k-1} e^{-x/\theta}/(\Gamma(k)\theta^k) 1_{\{x\ge 0\}}$ is the $\Gamma$-distribution with parameters $k = 2/\epsilon^2 - 1$ and $\theta = \epsilon^2\kappa/2$. The following proposition summarises their findings.
	
    \begin{Prop}[\cite{BCS25}]\label{prop:Rosenzweig_MacArthur_asymptotics}
	Let $(x_1, x_2) \in M_+$.
	\begin{enumerate}[label=(\roman*)]
		\item If $\epsilon > \sqrt 2$ then $X^{(x_1, x_2)}_t$ converges to $(0, 0)$ almost surely.
		\item If $\epsilon < \sqrt 2$ and $\Lambda > 0$ then there exists a unique invariant distribution $\mu_+$ on $M_+$ and $X^{(x_1, x_2)}_t$ converges to $\mu_+$ in total variation.
        \item If $\epsilon < \sqrt 2$ and $\Lambda < 0$ then $x_2(t) \to 0$ almost surely and the empirical occupation measure
        $$
            \Pi^{(x_1,x_2)}_t = \frac1t\int_0^t \delta_{X^{(x_1, x_2)}_s} ds
        $$
        converges weakly to $\mu_0$, almost surely.
    \end{enumerate}
	\end{Prop}

    Theorem~\ref{thm:main} provides the asymptotic behaviour of the law at time $t$ in the case of extinction, proving Conjecture 2.4 of~\cite{BCS25}.
    
	\begin{Prop}\label{prop:cv_Rosenzweig_MacArthur}
		Let $(x_1, x_2) \in M_+$. If $\epsilon < \sqrt 2$ and $\Lambda < 0$ then $X^{(x_1, x_2)}_t$ converges to $\mu_0$ in distribution.
	\end{Prop}

Note that $M_0 = \{x_1 x_2 = 0\}$ carries both $\mu_0 = \gamma_{\epsilon, \kappa}(x)dx \otimes \delta_0$ and the non-attracting invariant distribution $\delta_{(0, 0)}$. Hence
$$
\lim_{t \to +\infty}\sup_{(x_1, x_2) \in K} \FM{\delta_{(x_1, x_2)}P_t, \mu_0} = 0,
$$
as would be imposed by Assumption~\ref{asm:cv_to_mu0_single_measure} of Theorem~\ref{thm:single_measure}, fails when $(0, 0) \in K$.

\begin{Lem}\label{lem:ad_hoc_Lyapunov}
Set $M_{++} = \{x_1 > 0\}$. If $\epsilon< \sqrt 2$ then $W : M_{++} \to \mathbb R_+$ given by
 $$
 W(x_1,x_2)=e^{q(x_1+x_2)}x_1^{-p}
 $$
 where $0<q<2/(\kappa\epsilon^2)$ and $0<p<\min(\alpha q,2/\epsilon^2 -1)$ is a Lyapunov function for $(X_t)_{t \ge 0}$ restricted to $M_{++}$.
\end{Lem}

\begin{proof}
One has
\[
\begin{aligned}
	\frac{\mathcal L_eW}{W}(x_1,x_2)
	={}&
	-\left(\frac{q}{\kappa}-\frac{\epsilon^2q^2}{2}\right)x_1^2
	+\left(q+\frac{p}{\kappa}-\epsilon^2pq\right)x_1 \\
	&-\left(\alpha q-\frac{p}{1+x_1}\right)x_2
	-p\left(1-\frac{\epsilon^2(1+p)}{2}\right)
\end{aligned}
\]
and hence
\begin{align}\label{eq:LW/W_estimate}
\frac{\mathcal L_eW}{W}(x_1,x_2)
\le -A x_1^2
+\left(q+\frac{p}{\kappa}\right)x_1
-B x_2-D
\end{align}
where
\[
A=\frac{q}{\kappa}-\frac{\epsilon^2q^2}{2}>0,
\qquad B=\alpha q-p>0,
\qquad D=p\left(1-\frac{\epsilon^2(1+p)}{2}\right)>0.
\]
The right-hand side of~\eqref{eq:LW/W_estimate} is at most $-D/2$ under any of the following three conditions
\begin{align*}
\text{$x_1$ is sufficiently small,} \qquad\qquad \text{$x_1$ is sufficiently large,} \qquad\qquad x_2 \ge x_1 (q + p/ \kappa) / B.
\end{align*}
Thus $\mathcal L_e W \le -(D/2) W$ outside a compact set and the result follows.
\end{proof}

\begin{proof}[Proof of Proposition~\ref{prop:cv_Rosenzweig_MacArthur}]
It suffices to verify the assumptions of Theorem~\ref{thm:practical} with $n = 1$. Applying Proposition~\ref{prop:decay} to the Lyapunov function $W : M_{++} \to \mathbb R_+$ from Lemma~\ref{lem:ad_hoc_Lyapunov} yields~\ref{asm:tightness_practical}.	Assumption~\ref{asm:cv_to_M0_practical} holds by Proposition~\ref{prop:Rosenzweig_MacArthur_asymptotics}~(iii). Finally, the restriction of~\eqref{eq:Rosenzweig_MacArthur_SDE} to $M_0 \cap M_{++}$ is elliptic. Hence every point in $M_0 \cap M_{++}$ is accessible from $M_0 \cap M_{++}$~\cite[Prop.~6.33]{benaim22course} and a Doeblin point for $(X_t)_{t\ge 0}$ restricted to $M_0 \cap M_{++}$ by~\cite[Thm.~6.37]{benaim22course}, yielding Assumption~\ref{asm:cv_to_mu0_practical}.
\end{proof}

\subsection{Two-dimensional stochastic Kolmogorov system}\label{sec:Kolmogorov}

Finally, we study a general 2D stochastic Kolmogorov system, as in Example~2.4 of~\cite{HN18}. A distinctive feature of this model is that its extinction set can carry two attracting invariant distributions, causing $X^x_t$ to converge to the initial-state-dependent mixture $\mathbb P(\mathcal A_1^x) \mu_0^1 + \mathbb P(\mathcal A_2^x) \mu_0^2$. Consider the following SDE taking its values in $M = \mathbb R_+^2$
\begin{align}\label{eq:2d_kolmogorov}
    d x_1 &= x_1f_1(x_1, x_2)dt + x_1g_1(x_1, x_2)dE^1_t, & d x_2 &= x_2f_2(x_1, x_2)dt + x_2g_2(x_1, x_2)dE^2_t, 
\end{align}
where $f_1, f_2, g_1, g_2 : \mathbb R_+^2 \to \mathbb R$ are locally Lipschitz and the processes $(E^1_t)_{t\ge0}, (E^2_t)_{t\ge0}$ are given by $(E^1_t, E^2_t)^T = \Gamma (B^1_t,B^2_t)^T$. Here $(B^1_t)_{t\ge0}, (B^2_t)_{t\ge0}$ are independent Brownian motions and $\Gamma$ is a $2 \times 2$ real matrix. We denote the solution of~\eqref{eq:2d_kolmogorov} started from $(x_1, x_2) \in M$ by $X^{(x_1, x_2)}_t = (x_1(t), x_2(t))$. We work under Assumptions 1.1 and 1.4 of Hening and Nguyen~\cite{HN18}. Our goal in this section is not maximal generality but rather to illustrate how to use the existing theory of stochastic persistence/extinction to verify the assumptions of Theorem~\ref{thm:main}. We therefore also work under assumptions of~\cite[Prop.~3.2]{Ben18} and~\cite[Thm.~5.1]{Ben18} on the space $M$ as well as on the subspaces $\{ x_1 = 0\}$ and $\{x_2 = 0\}$. These assumptions are in particular verified by the stochastic Lotka–Volterra competitive model, see Example 2.1 of~\cite{HN18}. The coordinates $x_1$ and $x_2$ then represent two species' population levels and the persistence set $M_+ = \{x_1, x_2 > 0\}$ corresponds to the survival of both species, whereas the extinction set $M_0 = M \setminus M_+$ means that at least one species dies out.

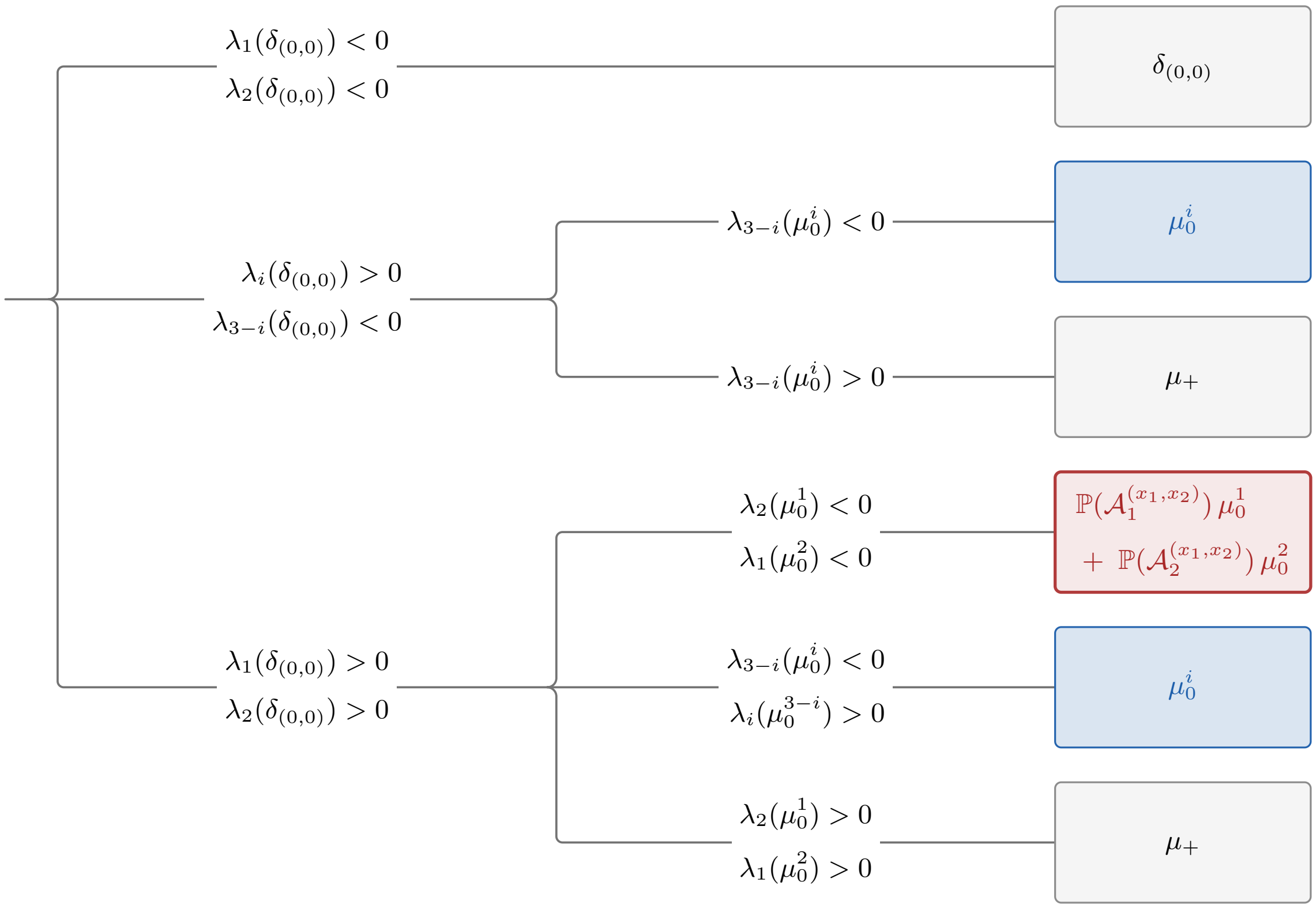
\begin{figure}[h!]
	\centering
	\noindent\makebox[\textwidth][c]{
		\resizebox{\textwidth}{!}{\begingroup%
\definecolor{boundarylimit}{RGB}{26,92,170}%
\definecolor{mixturelimit}{RGB}{168,38,38}%
\def\mixtureformula{\begin{aligned}%
    &\mathbb P(\mathcal A_1^{(x_1,x_2)})\,\mu_0^1\\[2pt]%
    &+\;\mathbb P(\mathcal A_2^{(x_1,x_2)})\,\mu_0^2%
\end{aligned}}%
\setbox0=\hbox{\footnotesize $\mixtureformula$}%
\edef\leaftextwidth{\the\wd0}%
\edef\leaftextheight{\the\ht0}%
\edef\leaftextdepth{\the\dp0}%
\pgfmathsetlengthmacro{\leafwidth}{1.1*(\leaftextwidth+6pt)}%
\pgfmathsetlengthmacro{\leafheight}{1.1*(\leaftextheight+\leaftextdepth+5pt)}%
\pgfmathsetlengthmacro{\leafxpadding}{(\leafwidth-\leaftextwidth)/2}%
\pgfmathsetlengthmacro{\leafypadding}{(\leafheight-\leaftextheight-\leaftextdepth)/2}%
%
\pgfmathsetlengthmacro{\leafleft}{\linewidth-\leafwidth}%
\pgfmathsetlengthmacro{\branchcolumn}{(\leafleft+6mm)/2}%
\pgfmathsetlengthmacro{\origincolumn}{(6mm+\branchcolumn)/2}%
\pgfmathsetlengthmacro{\invasioncolumn}{(\branchcolumn+\leafleft)/2}%
\pgfmathsetlengthmacro{\rowstep}{\leafheight+4mm}%
\pgfmathsetlengthmacro{\rowA}{-\rowstep}%
\pgfmathsetlengthmacro{\rowB}{\rowA-\rowstep}%
\pgfmathsetlengthmacro{\rowC}{\rowB-\rowstep}%
\pgfmathsetlengthmacro{\rowD}{\rowC-\rowstep}%
\pgfmathsetlengthmacro{\rowE}{\rowD-\rowstep}%
\pgfmathsetlengthmacro{\middleA}{(\rowA+\rowB)/2}%
\begin{tikzpicture}[
    branch/.style={draw=black!55, line width=0.7pt,
        rounded corners=2pt, line join=round, line cap=round},
    leaf/.style={draw=black!45, fill=black!4, rounded corners=2pt,
        line width=0.6pt, inner xsep=\leafxpadding, inner ysep=\leafypadding,
        minimum width=0pt, minimum height=0pt,
        text width=\leaftextwidth,
        text height=\leaftextheight, text depth=\leaftextdepth,
        anchor=east, align=center, font=\footnotesize, outer sep=0pt},
    single/.style={leaf, draw=boundarylimit!95,
        fill=boundarylimit!16, text=boundarylimit},
    mixture/.style={leaf, draw=mixturelimit!90,
        fill=mixturelimit!10, text=mixturelimit, line width=1pt},
    cond/.style={anchor=center, fill=white, inner xsep=2.5pt,
        inner ysep=2pt, outer sep=0pt, font=\footnotesize}
]
    \coordinate (root)     at (0pt,\middleA);
    \coordinate (forkRoot) at (6mm,\middleA);
    \coordinate (forkA)    at (\branchcolumn,\middleA);
    \coordinate (forkB)    at (\branchcolumn,\rowD);
    \coordinate (forkBtip) at ([xshift=-3pt]forkB);

    \node[leaf]    (delta) at (\linewidth,0pt)   {$\delta_{(0,0)}$};
    \node[single]  (leafA) at (\linewidth,\rowA) {$\mu_0^i$};
    \node[leaf]    (leafB) at (\linewidth,\rowB) {$\mu_+$};
    \node[mixture] (leafC) at (\linewidth,\rowC) {$\mixtureformula$};
    \node[single]  (leafD) at (\linewidth,\rowD) {$\mu_0^i$};
    \node[leaf]    (leafE) at (\linewidth,\rowE) {$\mu_+$};

    \node[cond] at (\origincolumn,0pt) {
        $\begin{aligned}
            \lambda_1(\delta_{(0,0)})&<0\\[2pt]
            \lambda_2(\delta_{(0,0)})&<0
        \end{aligned}$};
    \node[cond] at (\origincolumn,\middleA) {
        $\begin{aligned}
            \lambda_i(\delta_{(0,0)})&>0\\[2pt]
            \lambda_{3-i}(\delta_{(0,0)})&<0
        \end{aligned}$};
    \node[cond] at (\origincolumn,\rowD) {
        $\begin{aligned}
            \lambda_1(\delta_{(0,0)})&>0\\[2pt]
            \lambda_2(\delta_{(0,0)})&>0
        \end{aligned}$};

    \node[cond] at (\invasioncolumn,\rowA) {$\lambda_{3-i}(\mu_0^i)<0$};
    \node[cond] at (\invasioncolumn,\rowB) {$\lambda_{3-i}(\mu_0^i)>0$};
    \node[cond] at (\invasioncolumn,\rowC) {
        $\begin{aligned}
            \lambda_2(\mu_0^1)&<0\\[2pt]
            \lambda_1(\mu_0^2)&<0
        \end{aligned}$};
    \node[cond] at (\invasioncolumn,\rowD) {
        $\begin{aligned}
            \lambda_{3-i}(\mu_0^i)&<0\\[2pt]
            \lambda_i(\mu_0^{3-i})&>0
        \end{aligned}$};
    \node[cond] at (\invasioncolumn,\rowE) {
        $\begin{aligned}
            \lambda_2(\mu_0^1)&>0\\[2pt]
            \lambda_1(\mu_0^2)&>0
        \end{aligned}$};

    \begin{scope}[on background layer]
        \foreach \fork/\upper/\lower in {
            forkRoot/delta.west/forkBtip,
            forkA/leafA.west/leafB.west,
            forkB/leafC.west/leafE.west%
        } {
            \coordinate (braceTop) at (\fork |- \upper);
            \coordinate (braceBottom) at (\fork |- \lower);
            \draw[branch, sharp corners]
                (\upper)
                -- ([xshift=2pt]braceTop)
                arc[start angle=90, end angle=180, radius=2pt]
                -- ([yshift=3pt]\fork)
                .. controls ([yshift=1.34315pt]\fork)
                    and ([xshift=-1.34315pt]\fork)
                .. ([xshift=-3pt]\fork)
                .. controls ([xshift=-1.34315pt]\fork)
                    and ([yshift=-1.34315pt]\fork)
                .. ([yshift=-3pt]\fork)
                -- ([yshift=2pt]braceBottom)
                arc[start angle=180, end angle=270, radius=2pt]
                -- (\lower);
        }
        \draw[branch, sharp corners]
            (root) -- ([xshift=-3pt]forkA);
        \draw[branch, sharp corners]
            (forkBtip) -- (leafD.west);
    \end{scope}
\end{tikzpicture}%
\endgroup%
}
	}\par
	\caption{Decision tree to determine the limit in distribution of $X^{(x_1, x_2)}_t$ when $x_1, x_2 > 0$}
	\label{fig:decision_tree}
\end{figure}

The following quantities, defined for any invariant distribution $\mu$, determine the asymptotic behaviour of $(X_t)_{t\ge 0}$
$$
    \lambda_i(\mu) = \int_M \left(f_i(x_1, x_2) - \frac{\sigma_{ii} g_i^2(x_1, x_2)}{2} \right) d\mu(x_1, x_2) \hspace{1cm} (i = 1, 2)
$$
where $(\sigma_{ij})_{1\le i,j\le2}  = \Sigma= \Gamma\Gamma^T$. In particular, one has $\lambda_i(\delta_{(0,0)}) = f_i(0) - \sigma_{ii}g_i^2(0, 0)/2$ and if $\lambda_i(\delta_{(0,0)}) > 0$ then there exists an invariant distribution $\mu_0^i$ with $\mu_0^i(\{x_i > 0, x_{3-i} = 0\}) = 1$. The ellipticity condition contained in Assumption 1.1 of~\cite{HN18} ensures that $M_+$ carries at most one invariant distribution $\mu_+$. Figure~\ref{fig:decision_tree} is a decision tree yielding the limit in distribution of $X^{(x_1, x_2)}_t$ for $x_1, x_2 > 0$. Black leaves were established by Hening and Nguyen~\cite{HN18}, whereas coloured leaves are new and correspond to Proposition~\ref{prop:2d_kolmogorov}. The red leaf highlights a particularly striking regime: two attracting invariant distributions coexist on the extinction set and the limiting law of $X_t^x$ is a mixture whose weights retain a memory of the initial state.

\begin{Prop}\label{prop:2d_kolmogorov}
Let $x_1, x_2 > 0$.
	\begin{enumerate}[label=(\roman*)]
		\item If $\lambda_i(\delta_{(0,0)}) > 0, \lambda_{3-i}(\delta_{(0,0)})<0$ and $\lambda_{3-i}(\mu_0^i) < 0$, then
		$$
		X^{(x_1,x_2)}_t \text{ converges to } \mu_0^i \text{ in distribution}.
		$$
		\item If $\lambda_1(\delta_{(0,0)}), \lambda_2(\delta_{(0,0)}) > 0$ and $\lambda_2(\mu_0^1), \lambda_1(\mu_0^2) < 0$, then
		$$
		X^{(x_1,x_2)}_t \text{ converges to } \mathbb P(\mathcal A_1^{(x_1, x_2)})\mu_0^1 + \mathbb P(\mathcal A_2^{(x_1, x_2)})\mu_0^2 \text{ in distribution},
		$$
		where $\mathcal A_i^{(x_1, x_2)} = \left\{ \lim_{t \to +\infty} d(X^{(x_1, x_2)}_t, \{x_{3-i} = 0\}) = 0\right\}$.
		\item If $\lambda_1(\delta_{(0,0)}), \lambda_2(\delta_{(0,0)}) > 0$ and $\lambda_{3-i}(\mu_0^i) < 0, \lambda_i(\mu_0^{3-i}) > 0$, then
		$$
					X^{(x_1,x_2)}_t \text{ converges to } \mu_0^i \text{ in distribution}.
		$$
	\end{enumerate}
\end{Prop}

We follow the strategy laid out in Remark~\ref{rem:how_to_verify_assumptions}, i.e.~we apply the existing theory of stochastic persistence to verify the assumptions of Theorem~\ref{thm:practical}. In particular, we apply~\cite[Thm.~5.1]{Ben18} and~\cite[Prop.~8.3]{Ben18} to a carefully chosen persistence set $M_{++}$ to construct the Lyapunov function $W : M_{++} \to \mathbb R_+$ required by Assumption~\ref{asm:tightness_practical}.

\begin{proof}
We show how to obtain (i) from Theorem~\ref{thm:practical}. Assertions (ii) and (iii) follow from similar arguments. Fix $x_1, x_2 > 0$ and $i \in \{1, 2\}$ and assume
$$
\lambda_i(\delta_{(0,0)}) > 0, \lambda_{3-i}(\delta_{(0,0)})<0, \lambda_{3-i}(\mu_0^i) < 0.
$$
Applying~\cite[Thm.~5.1]{Ben18} and~\cite[Prop.~8.3]{Ben18} to the persistence set $M_{++} = \{ (y_1, y_2) \in M : y_i > 0 \}$ yields a Lyapunov function $W : M_{++} \to \mathbb R_+$ satisfying Assumption~\ref{asm:tightness_practical}. For this, observe that $\delta_{(0, 0)}$ is the only invariant distribution on $M_{00} = M \setminus M_{++}$ and hence the main condition ensuring persistence~\cite[Eq.~(33)]{Ben18} is equivalent to $\lambda_i(\delta_{(0,0)}) > 0$. Assumption~\ref{asm:cv_to_M0_practical} is established in Example 2.4 of~\cite{HN18}. Finally, the ellipticity assumptions under which we work yield Assumption~\ref{asm:cv_to_mu0_practical}.
\end{proof}

\section{Proofs}\label{sec:proof_main_result}

The proofs of the following three technical lemmas are relegated to the end of this section.

\begin{Lem}\label{lem:FM_mean} Let $(\Theta, m)$ be a probability space and let $\mu, \nu: \Theta \to \mathcal P(M)$ be measurable, i.e.~$\theta \mapsto \mu(\theta)(A)$ and $\theta \mapsto \nu(\theta)(A)$ are measurable for all measurable sets $A$. One has
	$$
	\FM{\int_\Theta \mu(\theta) dm(\theta), \int_\Theta \nu(\theta) dm(\theta)} \le \int_\Theta \FM{\mu(\theta), \nu(\theta)} dm(\theta),
	$$
	where $\FM{\cdot, \cdot}$ is the Fortet-Mourier distance of Definition~\ref{def:FM}.
\end{Lem}

\begin{Lem}\label{lem:uniform_in_time_Feller_continuity}
Under the standing Assumption~\ref{asm:CbM_Feller} one has
$$
\lim_{y \to x} \sup_{t \le T} \FM{\delta_y P_t, \delta_x P_t} = 0,
$$
for all $T > 0$ and $x \in M$.
\end{Lem}

\begin{Lem}\label{lem:separation_improvement}
   	If $C_1, C_2 \subset M$ are closed sets and $K \subset M$ is a compact set then $C_1 \cap C_2 \cap K = \emptyset$ implies that there exists $r > 0$ such that
   	$$
   	(C_1)_r \cap (C_2)_r \cap K = \emptyset
   	$$
where $(\cdot)_r = \{ x \in M : d(x, \cdot) < r\}$.
\end{Lem}

With these preliminary lemmas in hand, we turn to the proof of Theorems~\ref{thm:main} and~\ref{thm:practical}.

\begin{proof}[Proof of Theorem~\ref{thm:main}]
    	\setcounter{mainproofclaim}{0}
    	Let $\epsilon > 0$ be given. By Assumption~\ref{asm:tightness}, there exists $T > 0$ and a compact set $K_\epsilon$ such that
    	$$
    	\inf_{k \ge 1} \mathbb P\left(X^x_{kT} \in K_\epsilon \right) \ge 1 - \epsilon.
    	$$
    	By Lemma~\ref{lem:separation_improvement} and Assumption~\ref{asm:cv_to_M0} we have that there exists $r > 0$ such that $(M_0^i)_r \cap (M_0^j)_r \cap K_\epsilon = \emptyset$ for all $i \ne j$, where $(\cdot)_r = \{ x \in M : d(x, \cdot) < r\}$. By Assumption~\ref{asm:cv_to_mu0} there exists $T_0 > 0$ such that
    	$$
    	\sup_{x \in M_0^i \cap K_\epsilon} \FM{\delta_x P_t, \mu_0^i} \le \epsilon \text{ for all } t \ge T_0 \text{ and } i \in \{1, \ldots, n\}.
    	$$
    	By Lemma~\ref{lem:uniform_in_time_Feller_continuity}, for all $x \in M$ there exists a $0 < r_x < \min(r, \frac13 \min_{i\ne j} d(M_0^i \cap K_\epsilon, M_0^j \cap K_\epsilon) )$ such that
    	$$
    	\sup_{y \in B(x, r_x)} \sup_{t \le T + T_0} \FM{\delta_y P_t, \delta_x P_t} \le \epsilon.
    	$$
        Set $U_i = \bigcup_{x \in M_0^i \cap K_\epsilon} B(x, r_x)$ and $U = \bigcup_{i = 1}^n U_i$. We identify and prove four key claims before we can conclude our proof.
    	
    	\refstepcounter{mainproofclaim}\label{clm:main-convergence}
    	\textbf{Claim~\themainproofclaim:} One has $\lim_{k \to +\infty} d(X_{kT}, M_0) = 0$ almost surely where $M_0 = \bigcup_{i = 1}^n M_0^i$.
    	
    	\begin{proof}[Proof of Claim~\ref{clm:main-convergence}]
    	By Assumption~\ref{asm:cv_to_M0}, it suffices to show that $\mathbb P(\mathcal A_i^x \cap \mathcal A_j^x)=0$ for all $i \ne j$.
    	
        Let $\delta > 0$ be arbitrary but fixed. By Assumption~\ref{asm:cv_to_M0} and Lemma~\ref{lem:separation_improvement}, there exists $\tilde r > 0$ such that $(M_0^i)_{\tilde r} \cap (M_0^j)_{\tilde r} \cap K_\delta = \emptyset$. One has
    	$$
    	\lim_{k \to +\infty} 1_{\mathcal A_i^x \cap\mathcal A_j^x} 1_{\{ X^x_{kT} \in (M_0^i)_{\tilde r} \cap (M_0^j)_{\tilde r}\}} = 1_{\mathcal A_i^x \cap\mathcal A_j^x}
    	$$
    	so by dominated convergence
    	$$
    	\mathbb P(\mathcal A_i^x \cap \mathcal A_j^x) = \lim_{k \to +\infty} \mathbb P (\mathcal A_i^x \cap \mathcal A_j^x \cap \{ X^x_{kT} \in (M_0^i)_{\tilde r} \cap (M_0^j)_{\tilde r} \})
    	$$
    	and thus
    	\begin{align*}
    		\mathbb P(\mathcal A_i^x \cap \mathcal A_j^x) &\le \varlimsup_{k \to +\infty} \mathbb P (\mathcal A_i^x \cap \mathcal A_j^x \cap \{ X^x_{kT} \in (M_0^i)_{\tilde r} \cap (M_0^j)_{\tilde r} \cap K_\delta\}) + \varlimsup_{k \to +\infty} \mathbb P(X^x_{kT} \in K_\delta^c) \\
    		&\le 0 + \delta.
    	\end{align*}
    	We conclude by letting $\delta \to 0$.
    	\end{proof}
    	
    	\refstepcounter{mainproofclaim}\label{clm:main-outside-U}
    	\textbf{Claim~\themainproofclaim:} One has $\varlimsup_{k \to +\infty} \mathbb P_{x}(X_{kT} \in U^c) \le \epsilon$.
    	
    	\begin{proof}[Proof of Claim~\ref{clm:main-outside-U}]
    		The compact set $U^c \cap {K_\epsilon}$ is disjoint from the closed set $M_0$ as
    		$$
    		U^c \cap \underbrace{{K_\epsilon} \cap M_0}_{\subset U} \subset U^c \cap U = \emptyset,
    		$$
    		and thus $d(U^c \cap {K_\epsilon}, M_0) > 0$. So Claim~\ref{clm:main-convergence} implies $1_{U^c \cap {K_\epsilon}}(X^x_{kT}) \to 0$ almost surely
    		and, by dominated convergence, $\lim_{k \to \infty} \mathbb P(X^x_{kT} \in U^c \cap {K_\epsilon}) = 0$. We conclude by observing that $\mathbb P(X^x_{kT} \in U^c) \le \mathbb P(X^x_{kT} \in U^c \cap K_\epsilon) + \mathbb P(X^x_{kT} \in K_\epsilon^c)$ using Assumption~\ref{asm:tightness}.
    	\end{proof}
    	
    	 \refstepcounter{mainproofclaim}\label{clm:main-weights}
    	 \textbf{Claim~\themainproofclaim:} There exists $C_1(n) > 0$ such that
    	$$
    	\mathbb P(\mathcal A^x_i) - C_1(n)\epsilon < \varliminf_{k \to \infty} \mathbb P\left( X^x_{kT} \in U_i \right) \le \varlimsup_{k \to \infty} \mathbb P\left( X^x_{kT} \in U_i \right) < \mathbb P(\mathcal A^x_i) + C_1(n) \epsilon.
    	$$
    	
    	\begin{proof}[Proof of Claim~\ref{clm:main-weights}]
    		The union bound and $U_i \subset (M_0^i)_r$ imply
    		$$
    		\mathbb P(X^x_{kT} \in U_i) \le \mathbb P\left( \{ X^x_{kT} \in (M_0^i)_r \} \cap \mathcal A_i^x\right) + \sum_{j \ne i} \mathbb P\left( \{ X^x_{kT} \in (M_0^i)_r \} \cap \mathcal A_j^x\right).
    		$$
    		One has that $1_{\{X^x_{kT} \in (M_0^i)_r\} \cap \mathcal A_i^x}$ converges to $1_{\mathcal A_i^x}$ when $k \to +\infty$ and hence, by dominated convergence, one has $\mathbb P\left( \{ X^x_{kT} \in (M_0^i)_r \} \cap \mathcal A_i^x\right) \to \mathbb P\left(\mathcal A_i^x\right)$. For $j \ne i$ the same dominated convergence argument yields
    		$$
    		\varlimsup_{k \to \infty} \mathbb P\left( \{ X^x_{kT} \in (M_0^i)_r \} \cap \mathcal A_j^x\right) = \varlimsup_{k \to \infty} \mathbb P\left( \{ X^x_{kT} \in (M_0^i)_r \cap (M_0^j)_r \} \cap \mathcal A_j^x\right).
    		$$
    		and hence
    		\begin{align*}
    			\varlimsup_{k \to \infty} \mathbb P\left( \{ X^x_{kT} \in (M_0^i)_r \} \cap \mathcal A_j^x\right) &\le \varlimsup_{k \to \infty} \mathbb P\left( \{ X^x_{kT} \in (M_0^i)_r \cap (M_0^j)_r \cap K_\epsilon \} \cap \mathcal A_j^x\right) \\
    			&\qquad+ \varlimsup_{k \to +\infty} \mathbb P(X^x_{kT} \in K_\epsilon^c) \\
    			&\le \epsilon
    		\end{align*}
            using that $(M_0^i)_r \cap (M_0^j)_r \cap K_\epsilon = \emptyset$ so  $\varlimsup_{k \to \infty} \mathbb P\left( \{ X^x_{kT} \in (M_0^i)_r \cap (M_0^j)_r \cap K_\epsilon \} \cap \mathcal A_j^x\right) = 0$.
    		
             Hence $\varlimsup_{k \to +\infty} \mathbb P(X^x_{kT} \in U_i) \le \mathbb P(\mathcal A_i^x) + (n-1) \epsilon$.
    		
    		 For the other inequality, observe that
    		\begin{align*}
    			\varliminf_{k \to +\infty} \mathbb P(X^x_{kT} \in U_i) &\ge 1 - \varlimsup_{k \to +\infty} \mathbb P(X^x_{kT} \in U^c) - \sum_{j \ne i} \varlimsup_{k \to +\infty} \mathbb P(X^x_{kT} \in U_j) \\
    			&\ge 1 - \epsilon - (n - 1)^2 \epsilon - \sum_{j \ne i} \mathbb P(\mathcal A_j^x) \\
    			&= \mathbb P(\mathcal A_i^x) - ((n-1)^2+1)\epsilon
    		\end{align*}
    		using the fact that $\sum_{i=1}^n \mathbb P(\mathcal A_i^x) = 1$.
    	\end{proof}
    	
    	 \refstepcounter{mainproofclaim}\label{clm:main-mixture}
    	 \textbf{Claim~\themainproofclaim:} The measurable function $\nu: M \to \mathcal P(M)$ given by
    	$$
    	\nu(x) = \left\{
    	\begin{array}{cl}
    	\mu_0^i & \text{when } x \in U_i, \\
    	\mu_0^1 & \text{when }x \in U^c,
    	\end{array}
    	\right.
    	$$
satisfies
    	$$
    	\varlimsup_{k \to +\infty} \FM{\mathbb E[\nu(X^x_{kT})], \sum_{i = 1}^n \mathbb P(\mathcal A_i^x) \mu_0^i} \le C_2(n) \epsilon,
    	$$
    	for some $C_2(n) > 0$.
    	\begin{proof}[Proof of Claim~\ref{clm:main-mixture}]
    		Observe
    		$$
    		\mathbb E[\nu(X^x_{kT})] = \mathbb P(X^x_{kT} \in U^c) \mu_0^1+ \sum_{i = 1}^n \mathbb P(X^x_{kT} \in U_i) \mu_0^i.
    		$$
    		For measurable functions $f : M \to \mathbb R$ with $\Vert f\Vert_{\infty} \le 1$ one has
    		\begin{align*}
    			\mathbb E[\nu(X^x_{kT})](f) - \sum_{i = 1}^n \mathbb P(\mathcal A_i^x)\mu_0^i(f) &= \mathbb P(X^x_{kT} \in U^c) \mu_0^1(f) + \sum_{i = 1}^n \left(\mathbb P(X^x_{kT} \in U_i) - \mathbb P(\mathcal A_i^x)\right) \mu_0^i(f) \\
    			&\le \mathbb P(X^x_{kT} \in U^c) + \sum_{i = 1}^n \left|\mathbb P(X^x_{kT} \in U_i) - \mathbb P(\mathcal A_i^x)\right|
    		\end{align*}
    		Claim~\ref{clm:main-mixture} thus follows from Claims~\ref{clm:main-outside-U} and~\ref{clm:main-weights}.
    	\end{proof}
    	
    	 Having established these four intermediary claims, we are ready to conclude our proof of Theorem~\ref{thm:main}. Claims~\ref{clm:main-outside-U} and~\ref{clm:main-mixture} imply that there exists $N \ge 1$ such that for all $k \ge N$ one has
    	$$
    	\mathbb P(X_{kT} \in U^c) \le 2 \epsilon \text{ and } \FM{\mathbb E[\nu(X^x_{kT})], \sum_{i = 1}^n \mathbb P(\mathcal A_i^x) \mu_0^i} \le (C_2(n)+1)\epsilon.
    	$$
    	Let $t \ge NT + T_0$. There exist $k \ge N$ and $\Delta \in [T_0, T_0 + T)$ such that $t = kT + \Delta$. One has
    	\begin{align*}
    		\FM{\delta_x P_t, \sum_{i = 1}^n \mathbb P(\mathcal A_i^x) \mu_0^i} &\le \FM{\delta_x P_t, \mathbb E\left[ \nu(X^x_{kT})\right]} + \FM{\mathbb E\left[\nu(X^x_{kT})\right], \sum_{i = 1}^n \mathbb P(\mathcal A_i^x) \mu_0^i} \\
    		&\le \FM{\delta_x P_t, \mathbb E\left[ \nu(X^x_{kT})\right]} + (C_2(n)+1) \epsilon.
    	\end{align*}
    	Continuing our estimations, we have
    	\begin{align*}
    		\FM{\delta_x P_t, \mathbb E\left[ \nu(X^x_{kT})\right]} &= \FM{\mathbb E [\delta_{X^x_{kT}} P_\Delta], \mathbb E\left[ \nu(X^x_{kT})\right]} \\
    		&\le \mathbb E[\FM{\delta_{X^x_{kT}}P_\Delta, \nu(X^x_{kT})}] \\
    		&\le 2\mathbb P\left(X^x_{kT} \in U^c\right) + \sum_{i=1}^n \mathbb E[\FM{\delta_{X^x_{kT}}P_\Delta, \mu_0^i} 1_{\{X^x_{kT} \in U_i\}}].
    	\end{align*}
    	By construction, for all $y \in U_i$ there exists $x \in M_0^i \cap K_\epsilon$ with $\FM{\delta_y P_\Delta, \delta_x P_\Delta} \le \epsilon$. Furthermore, again by construction, one has $\FM{\delta_x P_\Delta, \mu_0^i} \le \epsilon$. Hence, by triangle inequality, one has $\FM{\delta_{y} P_\Delta, \mu_0^i} \le 2\epsilon$ for $y \in U_i$. Hence
    	$$
    	\FM{\delta_x P_t, \mathbb E\left[ \nu(X^x_{kT})\right]} \le 6\epsilon.
    	$$
Combining the previous estimates yields $\FM{\delta_x P_t, \sum_{i = 1}^n \mathbb P(\mathcal A_i^x) \mu_0^i} \le (C_2(n) + 7) \epsilon$, so letting $\epsilon \to 0$ concludes the proof.
\end{proof}

\begin{proof}[Proof of Theorem~\ref{thm:practical}]
Fix $x \in M_{++}$. It suffices to verify the assumptions of Theorem~\ref{thm:main}. By induction, it follows from~\ref{asm:tightness_practical} that
$$
P_{kT} W \le \rho^k W + \kappa \sum_{r =0}^{k-1} \rho^r,
$$
so by Markov's inequality
\begin{align*}
\mathbb P(W(X^x_{kT}) > a) \le \frac{\mathbb E[W(X^x_{kT})]}{a} \le \frac{W(x) + \kappa/(1 - \rho)}{a}.
\end{align*}
Assumption~\ref{asm:tightness} is thus satisfied by
$$
K_\delta = \{ y \in M_{++} : W(y) \le a(\delta) \} \text{ where } a(\delta) = (W(x) + \kappa / (1-\rho)) / \delta.
$$
Assumption~\ref{asm:cv_to_M0} follows from Assumption~\ref{asm:cv_to_M0_practical} and
$$
M_0^i \cap M_0^j \cap K_\delta \subset M_0^i \cap M_0^j \cap M_{++} = \emptyset \text{ for } i \ne j.
$$
Finally, applying~\cite[Thm.~8.15]{benaim22course} on $M_0^i \cap M_{++}$ yields $a_i, C_i > 0$ such that
$$
d_\mathrm{TV}(\delta_y P_t, \mu_0^i) \le C_i e^{-a_i t}(1 + W(y)) \text{ for all } y \in M_0^i \cap M_{++}.
$$
Because $M_0^i \cap K_\delta \subset M_0^i \cap M_{++}$ and $d_\mathrm{FM} \le d_\mathrm{TV}$ one has
$$
\sup_{y \in M_0^i \cap K_\delta} \FM{\delta_y P_t, \mu_0^i} \le C_i e^{-a_i t} \left(1 + \sup_{y \in M_0^i \cap K_\delta} W(y)\right).
$$
Hence Assumption~\ref{asm:cv_to_mu0} is satisfied.
\end{proof}

Finally, we provide the delayed proofs of this section's preliminary lemmas.
	
\begin{proof}[Proof of Lemma~\ref{lem:FM_mean}]
		Let $f$ be such that $ \Vert f \Vert_{\infty} \le 1 \text{ and } \Vert f \Vert_\mathrm{Lip} \le 1 $. One has
		\begin{align*}
			\FM{\int_\Theta \mu(\theta) dm(\theta), \int_\Theta \nu(\theta) dm(\theta)} &= \sup_{\Vert f \Vert_\infty, \Vert f \Vert_\mathrm{Lip} \le 1} \int_\Theta[\mu(\theta)(f)-\nu(\theta)(f)]\,dm(\theta) \\
			&\le \int_\Theta \sup_{\Vert f \Vert_\infty, \Vert f \Vert_\mathrm{Lip} \le 1} [\mu(\theta)(f)-\nu(\theta)(f)]\,dm(\theta) \\
			&= \int_\Theta \FM{\mu(\theta), \nu(\theta)} dm(\theta).
		\end{align*}
\end{proof}

\begin{proof}[Proof of Lemma~\ref{lem:uniform_in_time_Feller_continuity}] Assume by contradiction that there exists $\epsilon > 0$ such that
$$
\limsup_{y \to x} \sup_{t \le T} \FM{\delta_y P_t, \delta_x P_t} > \epsilon.
$$
Then for all $n \ge 1$, there exist $y_n \in M$ and $t_n \in [0, T]$ such that $d(y_n, x) < 1/n$ and $\FM{\delta_{y_n} P_{t_n}, \delta_x P_{t_n}} > \epsilon$. Up to extracting a subsequence, one can assume without loss of generality that $t_n \to t^*$. Because $d_\mathrm{FM}$ metrises convergence in distribution, Assumption~\ref{asm:CbM_Feller} implies
$$
\varlimsup_{n \to +\infty} \FM{\delta_{y_n} P_{t_n}, \delta_x P_{t_n}} \le \varlimsup_{n \to +\infty} \FM{\delta_{y_n} P_{t_n}, \delta_x P_{t^*}} + \varlimsup_{n \to +\infty} \FM{ \delta_x P_{t^*}, \delta_{x} P_{t_n}} = 0.
$$
This is a contradiction.
\end{proof}

\begin{proof}[Proof of Lemma~\ref{lem:separation_improvement}]
	Assume by contradiction that $(C_1)_r \cap (C_2)_r \cap K \ne \emptyset$ for all $r > 0$. Then we can construct a sequence $x_n \in K$ such that $d(x_n, C_1) < 1/n$ and $d(x_n, C_2) < 1/n$. Because $K$ is compact, up to extracting a converging subsequence, we can assume that $(x_n)_{n\ge 1}$ converges to a limit point $x$. By continuity of $d(\cdot, C_i)$ we have $d(x, C_i) = 0$. Closedness of $C_i$ implies $x \in C_i$ and hence $C_1 \cap C_2 \cap K \ne \emptyset$.
\end{proof}	

\paragraph{Acknowledgements.} M.B.~and L.H.~were funded by the Swiss National Science Foundation (SNSF) under grant no.~200020-219913. L.H. received further SNSF funding under grant no.~10009286.

\paragraph{AI declaration} GPT-6 Astra was used to improve the language and style of this article, review the literature and help with the LaTeX and TikZ code. It was also used as a proofreading tool and  provided the Lyapunov functions of Section~\ref{sec:applications}.  It was not used for any other mathematical contribution.

\bibliographystyle{alpha}
\bibliography{auto_update_biblio}
\end{document}